%% file: H-LU_pre_for_Max.tex
 \documentclass[graybox]{svmult}
 \usepackage{algorithm}
 \usepackage{romannum}
 
 \newcommand{\CD}{\mathcal{E}}
 \newcommand{\CE}{\mathcal{E}}
 \newcommand{\CH}{\mathcal{H}}
 \newcommand{\BHH}{\boldsymbol{H}}
 
 \newcommand{\BE}{\boldsymbol{E}}
 \newcommand{\BH}{\boldsymbol{H}}
 \newcommand{\BL}{\boldsymbol{L}}
 
 \newcommand{\skp}[1]{\left< #1 \right>}
 \newcommand{\norm}[1]{\left\|#1\right\|}
 \newcommand{\T}{\mathcal{T}}
 \newcommand{\CX}{\mathcal{Y}}

 \newcommand{\Bn}{\textbf{n}}
 \newcommand{\BB}{\textbf{B}}
 \newcommand{\BA}{\textbf{A}}
 \newcommand{\BU}{\textbf{U}}
 \newcommand{\BS}{\textbf{S}}
 \newcommand{\BV}{\textbf{V}}

 \usepackage{type1cm}        
 
 \newcommand{\blue}[1]{\textcolor{black}{#1}}
 \newcommand{\new}[1]{\textcolor{black}{#1}}
 \newcommand{\red}[1]{\textcolor{black}{#1}}
 
 \usepackage{makeidx}         
 \usepackage{graphicx}        
 \usepackage{multicol}        
 \usepackage[bottom]{footmisc}
 \newcommand{\RNum}[1]{\uppercase\expandafter{\romannumeral #1\relax}}

 \usepackage{newtxtext}       %
 \usepackage{newtxmath}       
 
 \makeindex             

 \usepackage{xcolor}

 \newcommand{\R}{\mathbb R}
 
 \usepackage{soul} 
 \usepackage{dsfont} 
 \input{commands}

\begin{document}
 	
 	\title*{Hierarchical LU preconditioning for the time-harmonic  Maxwell equations}
 	\author{Maryam Parvizi, Amirreza Khodadadian, Sven Beuchler, Thomas Wick}
 	\institute{Maryam Parvizi, Amirreza Khodadadian, Sven Beuchler, Thomas Wick \at Leibniz University Hannover, Institute of Applied Mathematics, Welfengarten 1, 30167 Hannover, 
 		Germany \email{\{parvizi,khodadadian,beuchler,thomas.wick\}@ifam.uni-hannover.de}
 		\and and Cluster of Excellence PhoenixD (Photonics, Optics, and
 		Engineering - Innovation Across Disciplines), Leibniz Universit\"at Hannover, Germany
 	}
 	\maketitle
 	\abstract{The time-harmonic Maxwell equations are used to study the effect of electric and magnetic fields on each other. Although the linear systems resulting from solving this system using FEMs are sparse, direct solvers cannot reach the linear complexity. In fact, due to the indefinite system matrix, iterative solvers suffer from  slow convergence.
 		In this work, we study the effect of using the inverse of $\CH$-matrix approximations  of the Galerkin matrices arising from  N\'ed\'elec's edge FEM discretization to solve the linear system directly. We also investigate the    impact of applying an $\mathcal{H}-LU$ factorization  as a preconditioner  and we study the number of iterations to solve the linear system using iterative solvers.
 	}\
 	\section{Introduction}
 	\label{sec:1}
 	The Maxwell system describes the behaviour of  electromagnetic fields.
 	N\'ed\'elec's edge finite element method is an efficient discretization technique to solve this equation which involves the curl-curl operator numerically \cite{monk2003finite, hiptmair2002finite}. Although the matrix arising from the linear system is sparse,  direct solvers cannot solve the problem with linear complexity. In the presence of  quasi-uniform meshes, the memory requirement of $O(N^{4/3})$ and computational time of $O(N^2)$ are expected  where the  problem size is  $N$ \cite{liu1992multifrontal}.

 	Matrices with full rank often can be approximated using low-rank matrices; but, it is not always applicable. Thus, it is desirable to present   a block-wise  partitioning of the matrix and approximate   appropriately chosen (using an admissibility condition) blocks by their low-rank decompositions.
 	Hierarchical matrices ($\mathcal{H}$-matrices), \cite{hackbusch2015hierarchical}, are block wise low-rank matrices that allow us to represent dense matrices with data sparse approximations and the  logarithmic-linear storage complexity, i.e.,  $\mathcal{O} (N m \log (N))$ where  $m$ is a parameter that controls the accuracy of the approximation.
 	
 	%
 	Besides the storage complexity reduction, another application of the $\CH$-matrix approximations   is to use them as a preconditioner to solve the system directly, or to use them as preconditioner to reduce the number of iterations in    Krylov-based iterative solvers based on matrix-vector multiplication, e.g., GMRES \cite{saad2003iterative}. 
 	In the time-harmonic case, the system matrix may become indefinite \red{and ill-conditioned}, in particular for high frequency \red{cases}. \red{In this regime, the usual factorization methods such as incomplete LU (iLU) do not lead to reliable results and converge to the exact solution poorly.} 
 	Then it is very difficult to design Galerkin discretizations \cite{MeSau11} and efficient 
 	iterative solvers \cite{ErnstGander} (see also \cite{HENNEKING202185} for recent studies).  

 	In this paper, we study the effect of applying hierarchical $LU$ decompositions, i.e., ${\mathcal{H}}-LU$ decompositions, as a preconditioner 
 	to solve the Maxwell equations using iterative solvers. \new{The idea of using $\mathcal{H}$ matrices for the curl-curl operators and magnetostatic model problems was introduced in \cite{bebendorf2009parallel}. A preconditioner based on the $\mathcal{H}$ matrices was used in \cite{bebendorf2012hierarchical} and \cite{ostrowski2010cal} for solving the Maxwell equations in the low-frequency regime. In \cite{faustmann2022h}, the authors  showed that the inverse of the Galerkin matrix corresponding to the FEM discretization of time-harmonic Maxwell equations can be approximated by an $\mathcal{H}$-matrix and proved the exponential convergence in the maximum block-rank.} \new{This approximation can be used to prove the existence of $\mathcal{H}$-$LU$ factorization without frequency restriction.} 
 	
 	In this work, in addition to addressing the advantage of iterative hierarchical preconditioning, we exploit the influence of applying an inverse of  $\CH$-matrix approximations as a preconditioner to solve the linear systems directly.
 	Our numerical tests also include studies on the influence of the wave number. \new{As observed, for both low and high frequency materials using ${\mathcal{H}}-LU$ factorization will lead to a fast and accurate convergence of the iterative solvers.}
 	
 	The paper is organized as follows. In Section \ref{section2}, we briefly introduce  the Maxwell system and present the  N\'ed\'elec's finite element discretization. The $\mathcal{H}$-matrices and how to compute the inverse of an $\CH$-matrix approximation for the Galerkin system matrix are explained in Section \ref{section3}. We also present an algorithm to compute the $\CH-LU$ approximation of a matrix, and use it as a preconditioner to solve the system directly.  In Section \ref{section4}, numerical results are presented to  substantiate the efficiency of the hierarchical matrix as a direct and iterative solver. Finally, the conclusions are given in Section \ref{conclusions}. 
 	\section{The Maxwell equations}\label{section2}
 	Denoting $\hat{\CH}$ the magnetic field intensity, $\CD$ the electric field density for the domain $\Omega\subset \mathbb{R}^d\,(d=2,3)$, we have the time-harmonic  Maxwell equations as
 	\begin{subequations}
 		\begin{alignat}{5}
 			\label{time dependent Gauss's law for mag}
 			\nabla \cdot (\beta \hat{\CH} ) &=0\quad \hspace{1cm}&&\text{in} \,\,\Omega\times \mathcal{I},\\
 			\label{time dependent Gauss's law for elec}
 			\nabla \cdot (\alpha \CE) &= \rho\quad \hspace{1cm}&&\text{in} \,\,\Omega\times \mathcal{I},\\ 
 			\label{time dependent Ampere's law}
 			\left(\alpha \frac{\partial }{\partial t }+\blue{\chi} \right) \CE\blue{-} \nabla \times \hat{\CH} &= \mathcal{F} \quad \hspace{1cm}&&\text{in} \,\,\Omega\times \mathcal{I},\\
 			\label{time dependent Faraday's law}
 			\beta \frac{\partial }{\partial t }\hat{\CH}+ \nabla \times \CE &=0 \quad \hspace{1cm}&&\text{in} \,\,\Omega\times \mathcal{I},
 		\end{alignat}  
 	\end{subequations}
 	where $\mathcal{F} $ is the applied electrical force, $\rho$ is the charge density, and  $\mathcal{I}=(0,T]$ is the time interval. In addition,  $\alpha$ and $\beta$ are the dielectric and magnetic permeabilities, \blue{and $\blue{\chi}$ is the conductivity constant.}  Considering an arbitrary  frequency $\omega$, with respect to time, the electric and magnetic fields can be represented as follows
 	\begin{subequations}
 		\begin{alignat}{5}
 			\label{2a}
 			\CE(x,t)&= e ^{-i \omega t} \boldsymbol{E}(x),\\
 			\label{2b}
 			\hat{\CH}(x,t)&= e ^{-i \omega t} \boldsymbol{H}(x),\\
 			\mathcal{F}(x,t)&=e ^{-i \omega t}\boldsymbol{F}(x).
 		\end{alignat}  
 	\end{subequations}
 	Replacing \eqref{2a} and \eqref{2b} into \eqref{time dependent Ampere's law} and \eqref{time dependent Faraday's law}, we \blue{obtain}\\[-8mm]
 	\begin{subequations}
 		\label{Maxwell system}
 		\begin{alignat}{5}
 			\label{Faraday's  law}
 			- \nabla \times \boldsymbol{H}-i  \omega\zeta  \boldsymbol{E} &= \textbf{F}(x)\quad \hspace{1cm}&&\text{in} \,\,\Omega,\\
 			\label{Ampere law}
 			\nabla \times \BE-i \omega \beta \BH &= 0\quad \hspace{1cm}&&\text{in} \,\,\Omega,
 		\end{alignat}  
 	\end{subequations}
 	where $\zeta:=  \blue{\alpha} + {i \blue{\chi} }/{\omega}$. 
 	Also, we  \blue{consider} a perfect conduction boundary condition (we surround $\Omega$ by a perfect bounded material i.e., $E\times \boldsymbol{n}=0$). Therefore, we have the following second-order operator for  \eqref{Maxwell system} 
 	\begin{align} \label{Maxwell-The electric field}
 		{\boldsymbol {\mathcal L}} \BE := \nabla  \times (\beta ^{-1}\nabla  \times \BE)- \kappa \BE= \boldsymbol{J}_S \qquad \text{in} \,\,\Omega ,
 	\end{align}
 	where  $\kappa:=\omega ^2 \zeta $ and $\boldsymbol{J}_S:=- i \omega \boldsymbol{F}$.
 	\subsection{Discretization by edge elements}
 	To present a Galerkin  weak formulation for \eqref{Maxwell-The electric field}, \blue{ we denote by 	 $\BL^2(\Omega)$ as the space of vector field with three entries from $L^2 (\Omega)$, i.e.,
 		$$\BL^2 (\Omega):= \left\lbrace \BU=(U_1, U_2,U_3) \quad:\quad U_i \in L^2 (\Omega), \,\,i=1,2,3\right\rbrace, $$}
 	with $\skp{\cdot\, , \cdot}$ as the inner product on this space, and continue with
 	defining the following space
 	$$\displaystyle \BH(\operatorname*{curl},\Omega):= \left\lbrace \mathbf{U} \in \BL^2 (\Omega ) \; \colon \; \nabla \times  \mathbf{U}\in \BL^2 (\Omega )  \right\rbrace,$$
 	equipped with the norm 
 	$$\norm{\textbf{U}}^2_ {\BH (\operatorname*{curl},\Omega)}:= \norm{\textbf{U}}^2_{\BL^2(\Omega)}+ \norm{\nabla \times \textbf{U}}^2_{\BL^2(\Omega)}.$$
 	Considering  homogeneous Dirichlet boundary conditions, the space 
 	$\BH_0 (\operatorname*{curl},\Omega) \subset \BH(\operatorname*{curl},\Omega)$ 
 	is introduced as follows
 	$$\BH_0 (\operatorname*{curl},\Omega):=\lbrace \textbf{U} \in \BL^2 (\Omega ) \; \colon \; \nabla \times  \textbf{U}\in \BL^2 (\Omega ) , \;  \BU \times  \Bn   =0 \;  \text{on}  \; \Gamma \rbrace. $$
 	Then, the weak formulation for \eqref{Maxwell-The electric field} can be written as
 	: find  $\BE \in \BHH_0 (\operatorname*{curl},\Omega)$ 
 	\begin{align}
 		\label{Maxwell weak form}
 		a(\BE, \Phi):= \skp{\nabla \times \BE, \nabla \times \Phi}  - \kappa \skp{\BE, \Phi}=
 		\skp{\boldsymbol{J}_S, \Phi}\qquad \forall \Phi \in \BHH_0 (\operatorname*{curl},\Omega).
 	\end{align}
 	Here, we should note that  $\kappa$ is not an eigenvalue of the  operator $\nabla \times \nabla \times$  \blue{\cite[Corollary .~{4.19}]{monk2003finite}}, in addition we have $\kappa \ne 0$\blue{, and we set $\beta=1$}.
 	The existence of the  unique  solution for the variational formulation \eqref{Maxwell weak form} is proven in \cite[Thm.~5.2]{hiptmair2002finite}, and  the following {\sl a priori} estimate is obtained 
 	\begin{equation}
 		\label{eq:apriori}
 		\norm{\BE}_{\BHH(\operatorname{curl},\Omega)} \leq C^* \norm{\boldsymbol{J}_S}_{L^2(\Omega)},
 	\end{equation}
 	where $C^*$ depends on $\Omega$ as well as $\kappa$. 
 	
 	For the \blue{discretization}, we consider  quasi-uniform mesh simplices ${\mathcal T}=\{T_1,\dots,T_{N_{\mathcal{T}}}\}$ where $T_j \in \mathcal{T}$ are open elements and denote
 	$h := \max_{T_j\in \mathcal{T}}{\rm diam}(T_j)$.
 	We assume  $\T$ is  a Ciarlet-regular mesh, i.e., it does not contain any hanging nodes. 
 	We also assume 
 	there is $\gamma >0$ such  that ${\rm diam}(T_\ell) \le \gamma\,|T_\ell|^{1/3}$ for all $T_\ell\in\mathcal{T}$.
 	In order to present the Galerkin FEM for \eqref{Maxwell weak form},  we consider lowest order N\'ed\'elec's $\operatorname*{curl}$-conforming elements of the first kind, i.e., 
 	\begin{align*}
 		\BV_h &:=\lbrace \mathbf{v } _h \in \BH(\operatorname*{curl},\Omega)\; \colon \; \mathbf{v } _h |_T \in \mathcal{N} _0 (T)   \quad \forall T \in \mathcal{T} \rbrace, \\
 		\BV_{h,0} &:=\BV_h \cap \BH_0(\operatorname{curl},\Omega),
 	\end{align*}
 	\blue{	where for all   $T \in \T _h $, the lowest order local N\'ed\'elec space of the first kind $\mathcal N _0(T)$  is  defined as
 		\cite{monk2003finite}} 
 	\begin{align*}
 		\blue{	\mathcal{N} _o (T) = \{ \boldsymbol{a} + \boldsymbol{b} \times \mathbf{x} \; \colon \; \boldsymbol{a} , \boldsymbol{b} \in \R ^3, \quad \mathbf{x} \in T\}.}
 	\end{align*}
 	We denote  $\CX_{h}:= \{\Phi_1,\ldots,\Phi_N\}$ as a basis with  $N$ as the dimension of $\BV_{h,0}$. 
 	This basis is uniquely defined by the property $\sigma(\Phi_i,e_j)=\delta_{ij}$,
 	where $e_j$ denotes an interior edge of the mesh and $\sigma(p,e)$ is the line integral of the tangential component of $p$ along $e$.
 	Then, the Galerkin FEM for \eqref{Maxwell weak form} is given as: 
 	Find  
 	$\BE _h \in \BV_{h,0}$ such that 
 	\begin{align}
 		\label{Galerkin discretization}
 		a(\BE_h, \Phi_h)=
 		\skp{\boldsymbol{J}_S , \Phi_h }\qquad \forall \Phi _h \in \BV_{h,0}.
 	\end{align}
 	This is equivalent to solve the following system 
 	\begin{align}
 		\label{eq:stiffness-matrix}
 		\BA x=b \quad\textrm{where}	\quad \BA=[\BA_{ij}]_{i,j=1}^N \quad\textrm{with}\quad
 		\BA _ {ij}:= a(\Phi _i , \Phi _j), \qquad  \Phi _j , \,\Phi _i \in \CX _ {h},
 	\end{align}
 	\blue{and the right hand side vector $b$ is defined as $b_j:= \skp{\boldsymbol{J}_S , \Phi_j }$, $j \in \{1,2,\ldots,N\}$. }
 	\section{$\mathcal{H}$-matrices and $\mathcal{H}$-matrix arithmetic }\label{section3}
 	$\mathcal{H}$-matrices are defined based on a partition $P$ generated from a clustering algorithm  that allows us to determine which blocks can be approximated by low-rank matrices or are  small \cite{hackbusch2015hierarchical}.   
 	
 	Applying  $\mathcal{H}$-matrix approximations  allows us to   store large matrices in the  format of low-rank block-wise matrices which could lead to logarithmic-linear storage complexity provided that  a proper method is used   to define the hierarchical structure  that results in the final block-wise format of the matrix. 
 	In the following lemma from \cite{faustmann2022h},  it is shown  that the inverse of the Galerkin matrix $\BA$ 
 	\eqref{eq:stiffness-matrix}
 	can be approximated using an $\CH$-matrix, and proven that this approximation converges exponentially with respect to the maximum block rank to $\BA$.\\[2mm]
 	\textbf{Definition 1}  \textit{[$\mathcal{H}$-matrices]\label{H-matrices}~A matrix $\mathbf{B}_ \CH  \in \mathbb{C} ^ {N\times N}$ is called an $\CH$-matrix, if for every
 		admissible block $(\tau , \sigma)$, we have the following factorization 
 		$$\mathbf{B}_ {\CH}|_{\tau \times \sigma}=\mathbf{X}_{\tau \sigma}  \mathbf{Y}^ H_{\tau \sigma},$$
 		of rank $r$	where $\mathbf{X}_{\tau \sigma} \in \mathbb{C} ^ {\tau\times r}$ and  
 		$\mathbf{Y}_{\tau \sigma} \in \mathbb{C} ^ {\sigma \times r}$.}\\
 	In order to use  the inverse of the $\CH$-matrix approximation of  $\BA$ as a preconditioner, first, we need to  find an $\CH$-matrix approximation for $\BA$, then we obtain the inverse using the iterative method of Schulz \cite{hackbusch2000sparse}.
 	\begin{lemma}\cite{faustmann2022h}
 		Let  $\BA$ be the Galerkin matrix defined in \eqref{eq:stiffness-matrix}. 
 		Then, there exists an $\mathcal{H}$-matrix approximation $\mathbf{B}_{\mathcal{H}}$ 
 		with the maximum block rank $r$ \blue{(rank of all the blocks of $\mathbf{B}_{\mathcal{H}}$ is either smaller than or equal to $r$)} such that 
 		\begin{align*}
 			\left\|\mathbf{A}^{-1} -\mathbf{B}_{\mathcal{H}}
 			\right\|_2 \le \bar{C}
 			h ^{-1} e^{-c(r^{1/4}/\ln r)},
 		\end{align*}
 		where $\bar{C}$ and $c$  are constants depending only on material parameters, and the properties of $\Omega$. 
 	\end{lemma}
 	In the definition of $\CH$-matrices, the low-rank blocks are determined based on the concept of $\eta$-admissibility defined in \cite{hackbusch2015hierarchical}. 	In the following, the mathematical definition of an $\CH$-matrix is given.\\[2mm]

 	
 	
 	Although computing the inverse of the $\mathcal{H}$-matrix approximation of the Galerkin system matrix  leads to logarithmic-linear complexity,  the computational cost to solve the linear system directly is still too high. Thus, we need to use another alternative to reduce the numerical cost such as the $\CH-LU$ factorization.
 	In the following, we present an algorithm from \cite[Sec. 2.9]{MR2451321}, and use it as a preconditioner to solve the linear systems.
 	\begin{algorithm}[ht!]\label{hlu}
 		\vspace{0.2cm}
 		\Romannum{1}) Compute the $\mathcal{H}$-matrix approximation of $\BA$, i.e., $\BA_\mathcal{H}$\\[2mm]
 		\Romannum{2}) Compute the $\mathcal{H}$- matrix $LU$ decomposition of $\BA_\mathcal{H}$ as follows\\[-2mm]
 		
 		\qquad 1) \textbf{for} $j=1,\ldots,N$\\[-2mm]
 		
 		\qquad\quad  \textbf{for} $k=1,\ldots,j-1$\\[-2mm]
 		
 		\qquad\qquad Solve the system $\displaystyle{\sum_{i=1}^{k}}L_{ji}\,U_{ik}={(\BA_ \mathcal H)}_{jk}$ to get $L_{jk}.$\\
 		
 		\qquad 2) Compute $L_{jj}$ and $U_{jj}$ by~~~ $\displaystyle{L_{jj} U_{jj}=A_{jj}-\sum_{i=1}^{j-1} L_{ji} U_{ij}}$.\\
 		
 		\qquad 3) \textbf{for} $k=j+1,\ldots,N$\\[-2mm]
 		
 		\qquad\qquad Solve the system $\displaystyle{\sum_{i=1}^{j}}L_{ji}\,U_{ik}={(\BA_ \mathcal H)}_{jk}$ to get $U_{jk}.$\\
 		
 		\Romannum{3})  \textbf{for} $i=1,\ldots,\text{MaxIt}$\\[-2mm]
 		
 		\qquad Compute~~$r=b-\BA\,x$\qquad \text{and} \qquad \texttt{err}=$\norm{r}_2/ \norm{b}_2$
 		\\[-2mm]   	
 		
 		\qquad Compute~~$x=x+U^{-1} \left(L^{-1}\,r\right)$.\\[-2mm]
 		
 		\qquad \textbf{if}  $ \texttt{err}<\texttt{TOL}$ \quad $\textbf{break}$\\[-1mm]
 		
 		\Romannum{4}) Compute $\texttt{error}=\norm{\BA\,x-b}_2$.  	
 		\caption{ $\mathcal{H}-LU$ decomposition and  application in  preconditioning of a simple iterative  solver for \eqref{Galerkin discretization}.}
 	\end{algorithm}
 	\section{Main algorithm and numerical experiments}
 	\label{section4}
 	In this section, we first present the main algorithm of this work, and then we study three numerical examples to solve the linear system arising from the Maxwell equations using an  $\mathcal{H}$-matrix approximation as a preconditioner. For this, we employ the geometrically balanced cluster tree presented in  \cite{grasedyck2004adaptive}, and we set the admissibility parameter $\eta=2$.
 	We use a truncated singular value decomposition (SVD) with different ranks $r$ to compute  $\BB_{\mathcal{H}}$  as the inverse of $\BA_{\CH}$ obtained from the Schulz iteration. In other words, for admissible blocks $(\tau , \sigma)$, we have $\BA _ {\CH}|_ {\tau \times \sigma} := \BU _r\BS_r \BV^T_r $  where 
 	$\BU _r \in \R ^ {\tau \times r}$, $\BS_r \in  \R ^ {r \times r}$ and  $\BV _r\in \R ^{\sigma \times r}$ are the first $r$ columns of $\BU$, $\BS$ and  $\BV$, respectively. Then, we find the inverse $\BB_{\CH}$ for the matrix  $\BA _ {\CH}$.\\[-12mm]
 	\subsection{Example 1: a unit box}	
 	The geometry is $\Omega=(0,1)^3$ and $\boldsymbol{J}_S=[0,0,1]^\top$. The coarse mesh consists of 6 tetrahedrons. This mesh is uniformly refined $k$ times, $k=2,\ldots,7$. 
 	All computations are performed in MATLAB with 125 cores. 
 	Table \ref{tab:3} displays the iteration numbers of the preconditioned GMRES method with the described $\mathcal{H}$-matrix preconditioner for different values of $\kappa$. The GMRES method is stopped if a relative accuracy of $\blue{TOL}=10^{-5}$ of the residual is reached. In all experiments, the parameter $\beta=1$ is chosen.
 	After 100 iterations, we restart GMRES.
 	The results show that the $\mathcal{H}$-matrices can be used as an efficient preconditioner if $\kappa\leq 20^2$.
 	For higher frequencies, the iteration numbers grow in some, but not all levels $k$.
 	This is due to the computation of the $LU$ decomposition of the $\mathcal{H}$-matrix. The approximation of the original matrix by the $\mathcal{H}$-matrix is still very good, also in the case of $k=5$ and \blue{$\kappa=900$}.
 	
 	
 	\begin{table}[htbp]
 		\centering
 		\begin{tabular}{|c|c|cccccc|}
 			\hline
 			\blue{$\kappa~~$} & $N_{\text{dof}}$ &  $~~25~~$ & $~~100~~$ & $~~225~~$ & $~~400~~$ & $~~625~~$ & $~~900~~$ \\[0.5mm]
 			\hline 
 			\blue{$k=2$} & $98$ & $1$ & $1$ & $1$ & $1$ & $1$ & $1$ \\[2mm]
 			\blue{$k=3$} & $604$ & $2$ & $2$ & $2$ & $2$ & $2$ & $2$ \\[2mm]
 			\blue{$k=4$} & $3\,184$ & $2$ & $4$ & $5$ & $5$ & $4$ & $3$ \\[2mm]
 			\blue{$k=5$} & $41\,024$ & $3$ & $4$ & $8$ & $9$ & $11$ & $>3000$ \\[2mm]
 			\blue{$k=6$} & $238\,688$ & $3$ & $6$ & $20$ & $4$ & $48$ & $25$ \\[2mm]
 			\blue{$k=7$} & $1\,807\,264$ & $4$ & $5$ & $7$ & $19$ & $46$ & $93$ \\[2mm]
 			\hline
 		\end{tabular}
 		\caption{GMRES iterations numbers for Example 1.}
 		\label{tab:3}
 	\end{table}
 	\subsection{Example 2: two boxes}	
\begin{figure}[htbp]
 		\centering
 		\includegraphics[width=0.8\textwidth]{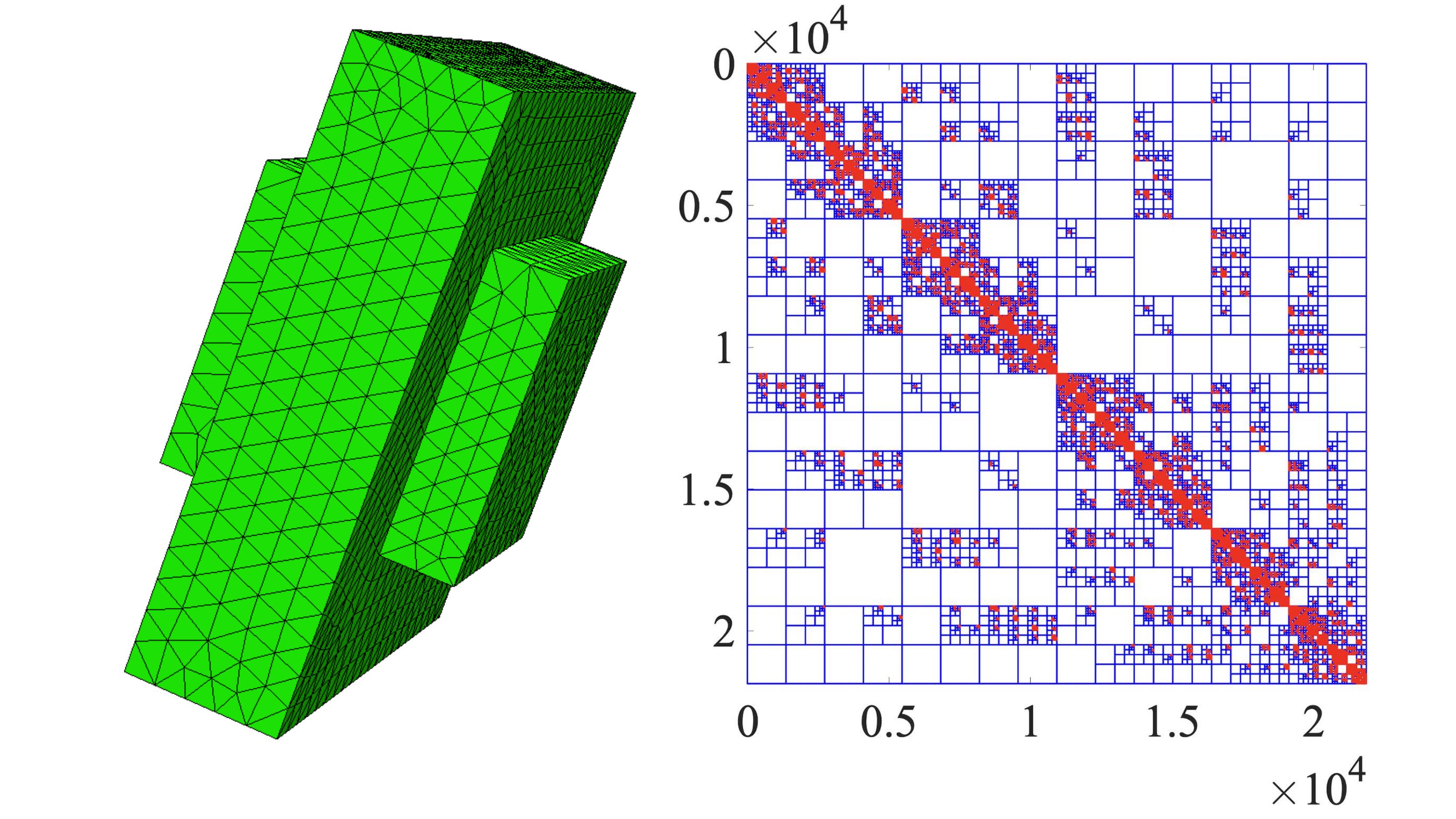}
 		\caption{Example 2. The computational geometry (left) and the inverse of $\CH$-matrix approximation of the stiffness matrix \eqref{Galerkin discretization} $\mathcal{H}$-matrix clustering (right) for $N_{\text{dof}}=22\,001$.}
 		\label{Example1:domain}
 	\end{figure}

\begin{figure}[htbp]
 		\centering
 		\includegraphics[width=0.8\textwidth]{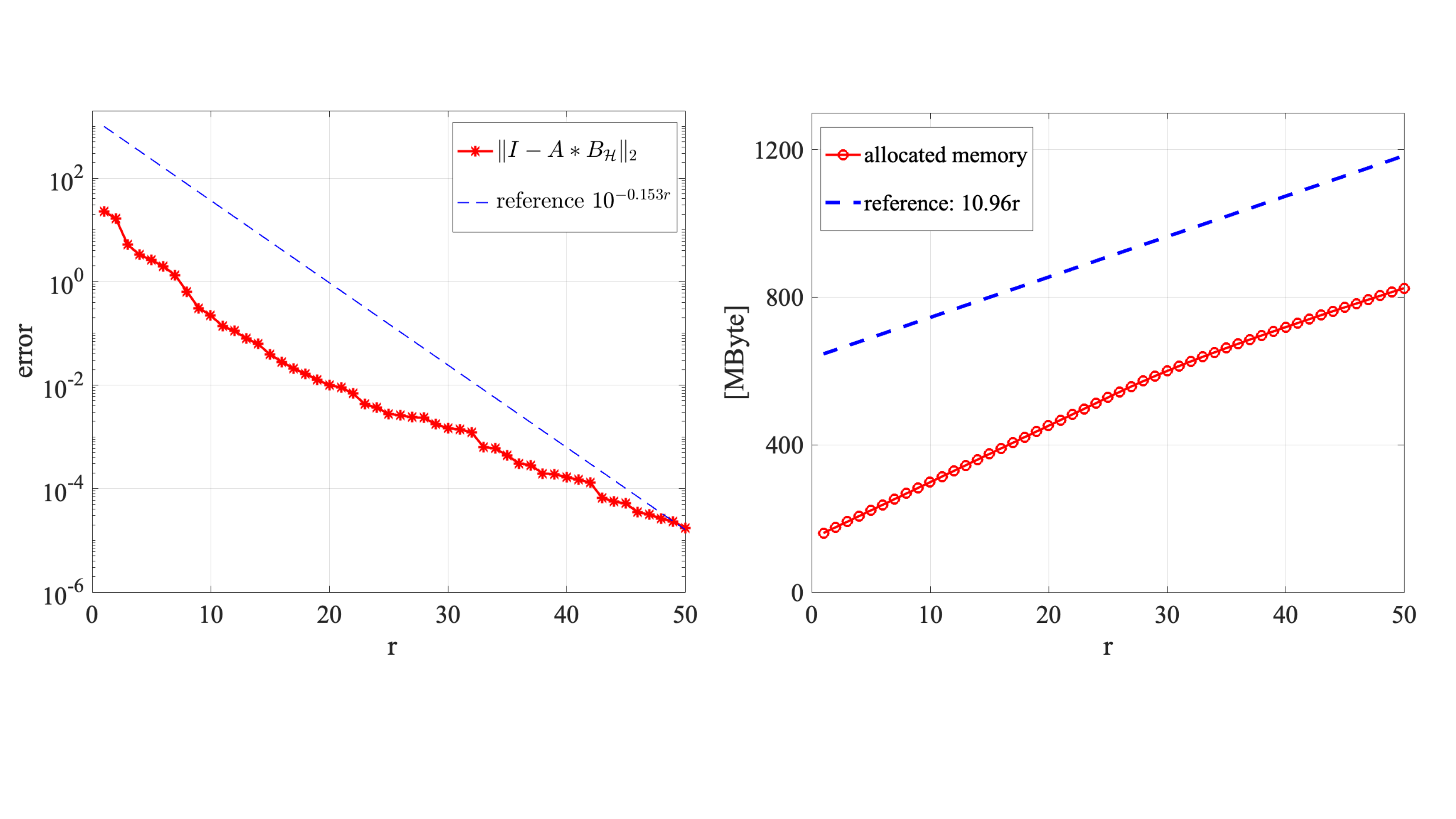}
 		\caption{Approximation error of $B_\mathcal{H}$ in Example 2 for a stiffness matrix with $N_{\text{dof}}=22\,001$ (left) and the relative allocated memory (right).}
 		\label{Example1:decay}
 	\end{figure}

 	Here the geometry consists of two boxes, i.e., $\Omega := (-1,1) \times (-1,1)\times [-1, -2) \cup (-2,2) \times [1,2) \times (-1,1) \cup (-2,2) \times [-1,1] \times [-1,1]\cup (-1,1) \times (-1,1) \times [1,2) \cup (-2,2) \times [-1,2) \times (-1,1)$.  We set the parameters    $\kappa=1$, $\boldsymbol{J}_S=[1,1,1]^\top$, 
 	and $\beta=1$.
 	The computational domain with and the inverse of  $\mathcal{H}$-matrix approximation of the stiffness matrix  is shown in Figure \ref{Example1:domain}. We have 24\,440 admissible blocks, 48\,404 small blocks and the  depth is 16.
 	The decay of the approximation error versus $r$ and the corresponding allocated memory are shown in Figure \ref{Example1:decay}. As shown, using higher $r$ gives rise to a reliable inverse of  the $\CH$-matrix approximation. The computed $\BB_{\mathcal{H}}$ can be used as a preconditioner  to solve the linear system \eqref{Galerkin discretization} directly.
 	
 	\subsection{Example 3: a magnet surrounded by air}	
\begin{figure}[htbp]
 		\centering
 		\includegraphics[width=0.8\textwidth]{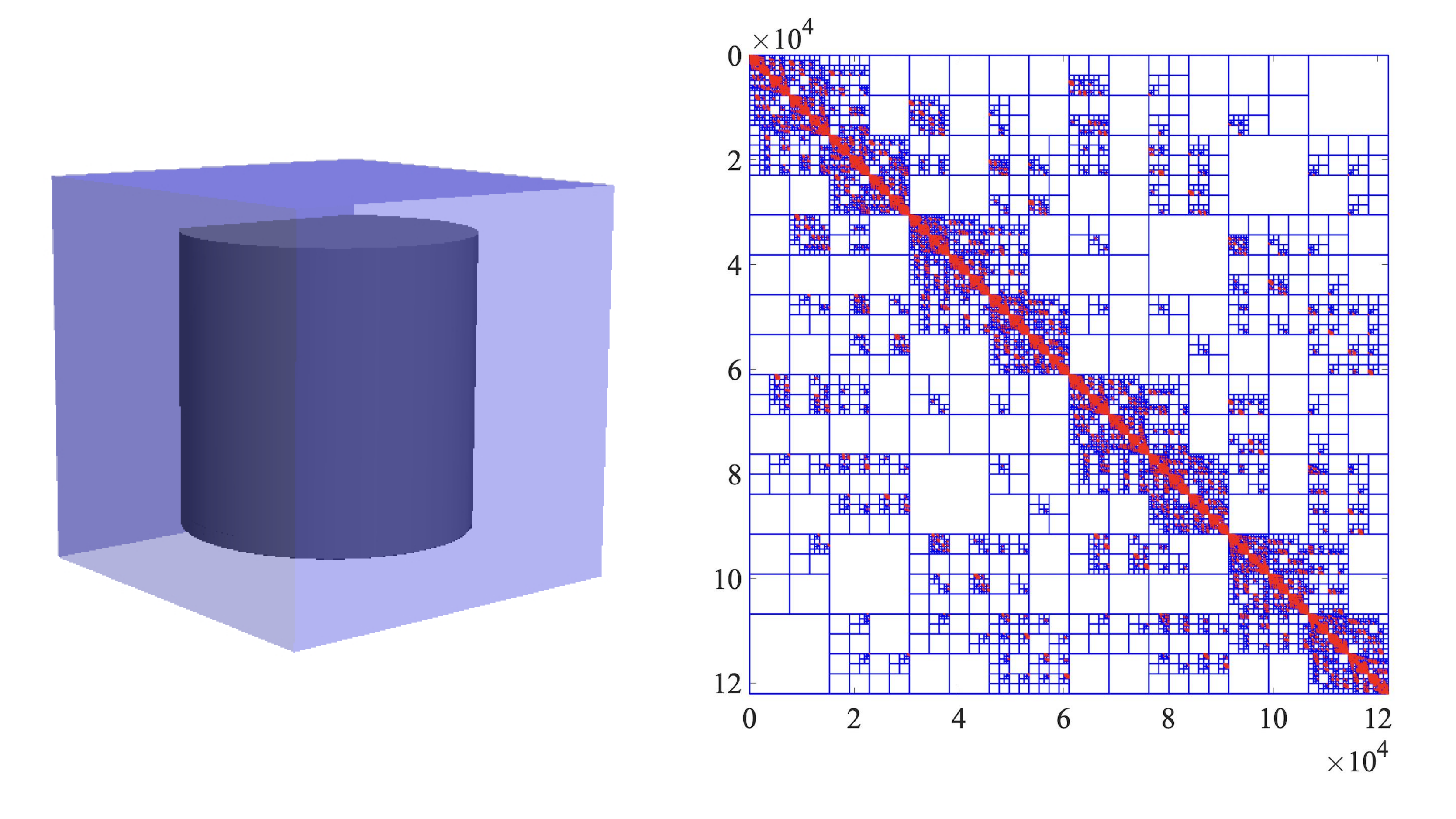}
 		\caption{Example 3. The computational geometry (left) and  $\mathcal{H}$-matrix approximation of $\BA$ for $N_{\text{dof}}=122\,202$ (right).}
 		\label{Example2}
 	\end{figure}

\begin{figure}[htbp]
 		\centering
 		\includegraphics[width=0.8\textwidth]{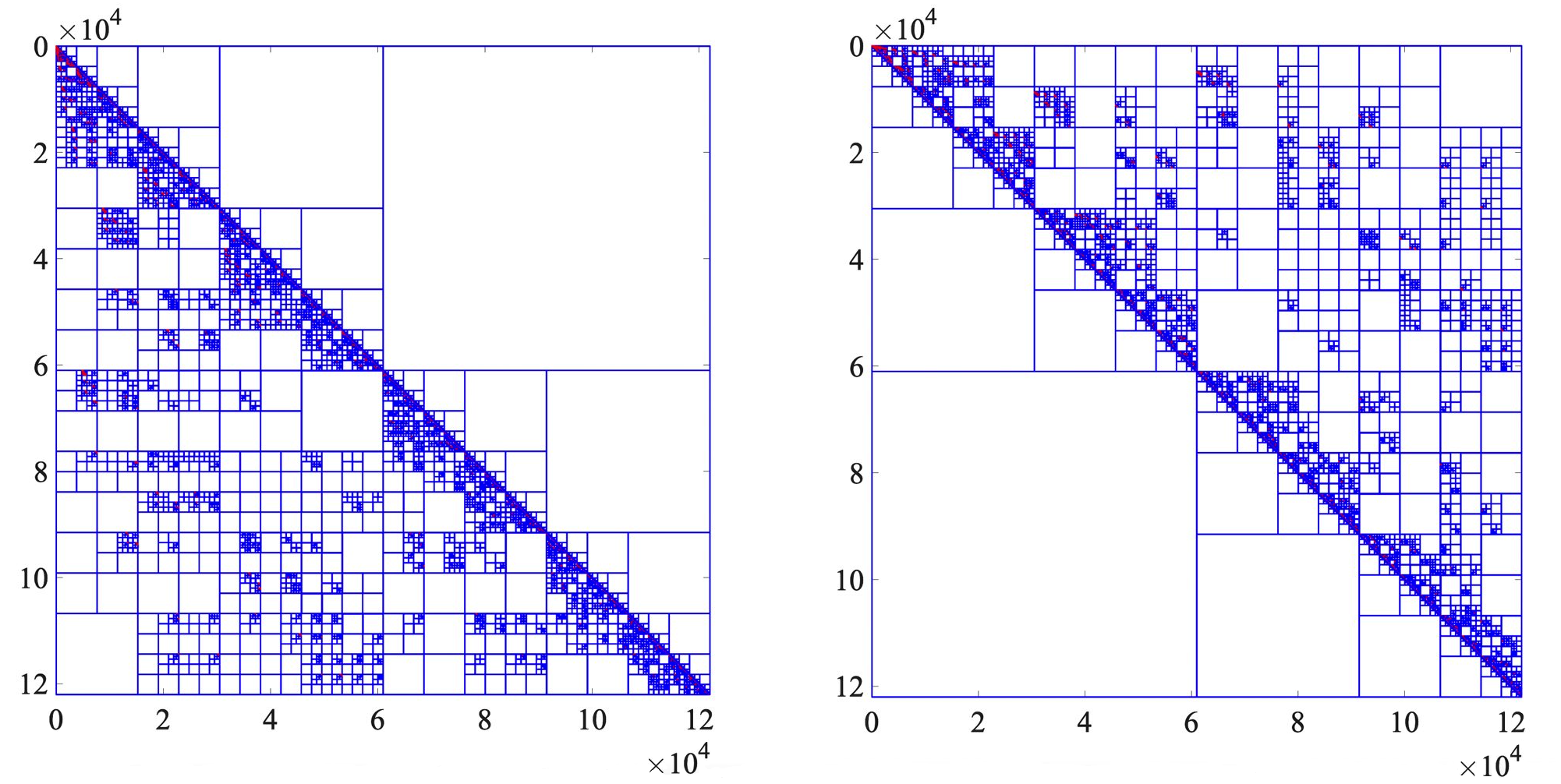}
 		\caption{The  $\CH$-LU decomposition for the stiffness matrix resulting from Example 3 corresponding to Figure \ref{Example2}.}
 		\label{Example2:HLU}
 	\end{figure}

 	We consider a magnet surrounded by air where the box is $1\times 1\times 1$ and the magnet dimension is $0.5\times 0.5\times 0.75$. We set $\kappa=10$, $\beta=10$ and $\boldsymbol{J}_S=[10,10,10]^\top$. 
 	The geometry and the  $\mathcal{H}$-approximation of $\BA$ for $N_{\text{dof}}=122\,202$ is shown in Figure \ref{Example2}. 
 	In this approximation, we have 215\,964 admissible blocks, 402\,451 small blocks, and the depth is 15. The ${\mathcal{H}}-$LU decomposition of $\BA$ is given in Figure \ref{Example2:HLU}. 
 	We use Algorithm 1 for different $N_{\text{dof}}$ to solve the linear system. Table \ref{comp2} shows the results for different matrices using $TOL=1\times 10^{-8}$. For all cases, the negligible $\texttt{error}$ indicates the accuracy the   method, and the elapsed CPU time points out its efficiency.  For the last two examples, we used \textit{Netgen/ngsolve} \cite{schoberl2017netgen} to produce the meshes (denoting different $N_\text{dof}$) and assembling the matrices of \eqref{Galerkin discretization}.
 	
 	\vspace{-5mm}
 	
 	\begin{table}[htbp]
 		\caption{The iterative solver preconditioned by $\mathcal{H}-LU$ used to solve the Maxwell equations for different $N_{\text{dof}}$.}
 		\vspace{1mm}
 		\centering
 		\begin{tabular}{|l| c c c c c c c c c c|}
 			\hline
 			$N_{\text{dof}}$ & ~~5\,492 ~~        &  ~~8\,050 ~~    &   ~~ 13\,602~~    &   ~~ 33\,933~~    &~~ 48\,811~~   &~~ 70\,133~~  & ~~ 78\,603~~  &~~ 96\,846~~  &~~ 129\,200~~  & ~~ 304\,309 \\\hline
 			\texttt{error}   & ~~7.22\,e\,-9~  & ~1.07\,e\,-9~  & ~5.80e-9~ &~2.23e-8~ & ~5.11e-8~& ~9.06e-8~& ~5.98e-10~&~4.38e-10~&~5.04e-10~&~5.78e-9~\\[4mm]
 			time [s]  & 11.52   & 35.24  & 60.52 &151.71 &  256.18 & 546.9& 602.73& 481&820.9&2335\\[4mm]
 			iterations &   3  &  3  &  3 & 3 &4  & 4& 5&5& 6&7\\[2mm]
 			\hline
 		\end{tabular}\label{comp2}
 		
 	\end{table}
 	
 	\section{Conclusion}
 	\label{conclusions}
 	In this work, the $\CH$-matrix approximation was used to solve the time-harmonic Maxwell equations. As a direct solver, the inverse of the hierarchical matrix approximation of the linear system was employed as a preconditioner, where an accurate approximation was achieved. Additionally,  we then employed an $\mathcal{H}-LU$ factorization as a preconditioner. 
 	
 	In both cases, the use of  $\CH$ matrix approximations could reduce the computational cost and increase the accuracy of the solution. \red{The $\mathcal{H}$ matrices  can be coupled with the domain decomposition to take advantage of both approaches, i.e., to reduce the complexity  and accelerate the convergence of the  iterative solvers. This possibility will be addressed in the future papers.}
 	\begin{acknowledgement}
 		The authors acknowledge  
 		the Deutsche Forschungsgemeinschaft (DFG) under Germany Excellence Strategy within  the Cluster of Excellence PhoenixD (EXC 2122, Project ID 390833453). Maryam Parvizi is funded by  Alexander von Humboldt Foundation project named $\mathcal{H}$-matrix approximability of the inverses for FEM, BEM and FEM-BEM coupling of the electromagnetic problems. Finally, the authors thank Sebastian Kinnewig for fruitful discussions. \\
 	\end{acknowledgement}
 	
 	\bibliographystyle{abbrv}
 	\bibliography{lit}
 	
 \end{document}

%% file: commands.tex
\definecolor{lightblue}{RGB}{0, 142, 255}

\definecolor{sblue}{RGB}{0, 10, 255}